\documentclass[11pt,twoside,leqno]{article}
\usepackage{amssymb}
\usepackage{epsfig,latexsym,color,colordvi}
\usepackage{mathrsfs}
\usepackage{latexsym}
\usepackage{eufrak}
\makeatletter
\@addtoreset{equation}{section}

\makeatother
\renewcommand\thefigure{\thesection.\@arabic\c@figure}
\renewcommand\thetable{\thesection.\@arabic\c@table}

\reversemarginpar

\hsize \textwidth      
\advance \hsize by -\marginparwidth
\evensidemargin \oddsidemargin 

\advance\hoffset by 5mm 

\newcommand{\mysection}[1]{\bigskip\par\refstepcounter{section}
{\bf\thesection.~#1}}
\newcounter{theorem}
\newenvironment{theorem}
{\bigskip\par\indent{\sc Theorem} \refstepcounter{theorem}\arabic{theorem}. \it}
{\bigskip\par}
\newcounter{remark}
\newenvironment{remark}
{\bigskip\par\indent{\sc Remark} \refstepcounter{remark}\arabic{remark}. }
{\refstepcounter{remark}\bigskip\par}
\newcounter{lemma}
\newenvironment{lemma}
{\bigskip\par\indent{\sc Lemma} \refstepcounter{lemma}\arabic{lemma}. \it}
{\refstepcounter{lemma}\bigskip\par}
\newcounter{claim}

\newcounter{definition}

\newcounter{corollary}

\newcounter{proposition}
\newenvironment{proposition}
{\bigskip\par\indent{\sc Proposition} \refstepcounter{proposition}\arabic{proposition}. \it}
{\refstepcounter{proposition}\bigskip\par}
\def\proof{{\sc Proof}. }

\def\mymonth{\ifcase\month\or January\or February\or March\or April
\or May\or June\or July\or August\or September
\or October\or November\or December\fi\ }

\def\cf{{\it cf.\/}\ }  
\def\st{\,\vert\,}
\def\qed{ {\small$\Box$}}
\def\eps{\varepsilon}
\def\R{{{\rm I}\!{\rm R}}}
\def\N{ { {\rm I}\!{\rm N}} }
\def\Z{{{\rm Z}\mkern-5.5mu{\rm Z}}}
\def\P{{{\rm I}\!{\rm P}}}
\def\E{{{\rm I}\!{\rm E}}}
\def\ind{{{\rm 1\mkern-1.5mu}\!{\rm I}}}
\def\w{\omega}

\def\card#1{\hbox{\rm card}\left(#1\right)}
\def\log{\hbox{\rm Log}\,}
\def\linf{\mathop{\underline{\lim}}}
\def\lsup{\mathop{\overline{\lim}}}
\def\ord{\kern0.1em o\kern-0.02em{}_{\ds\breve{}}\kern0.1em}
\def\norm#1{\left\Vert #1\right\Vert }
\def\abs#1{\left\vert\, #1\,\right\vert }
\def\ds{\displaystyle}

\def\vv{{\hbox{w}}}
\def\VV{{\hbox{W}}}
\def\svv{{\hbox{\scriptsize w}}}
\def\sVV{{\hbox{\scriptsize W}}}

\title{\vspace{-2Em}\vspace{-0.5ex}\normalsize\bf THE POINT OF VIEW OF THE
PARTICLE ON THE LAW OF LARGE NUMBERS FOR RANDOM WALKS IN A MIXING RANDOM ENVIRONMENT}
\author{{\sc By Firas Rassoul--Agha}\\
\\ {\it Courant Institute}
}
\date{}

\begin{document}

\maketitle \vspace{-1Em}
\hfill
\begin{minipage}{0.8\textwidth}
\footnotesize\hspace*{4ex}The point of view of the particle is an 
approach that has proven
very powerful in the study of many models of random motions in random
media. 
We provide a new use
of this approach to prove the law of large numbers in the case of one or 
higher-dimensional, finite range, transient random walks in mixing random 
environments.
One of the advantages of this method over what has been
used so far is that it is not restricted to i.i.d. environments.
\end{minipage}
\hfill
\renewcommand{\thefootnote}{}
\addtolength{\skip\footins}{1EM}
\footnotetext{Supported by NSF Grant DMS--010343.}

\bigskip

\bigskip

\bigskip

\mysection{Introduction.}
Originating from the physical sciences, the subject of random media has
gained much interest over the last three decades.
One of the fundamental models in the field is random walks
in a random environment. The main purpose of this work is to prove the
law of large numbers for a certain class of random walks in a mixing
random environment.
In this model, an environment is a collection of
transition probabilities
$\w=(\pi_{ij})_{i,j\in\Z^d}\in[0,1]^{\Z^d\times\Z^d}$, with
$\sum_{j\in\Z^d}\pi_{ij}=1$.
Let us denote by $\Omega$, the space of all such transition 
probabilities.
The space $\Omega$ is equipped with the canonical product
$\sigma$-field $\mathfrak{S}$, and with the natural shift
$(T^k\w)_{i,j}=\w_{k+i,k+j}$, for $k\in\Z^d$.
Here, $\w_{ij}$ stands for the $(i,j)^{\hbox{th}}$
coordinate of $\w\in\Omega$. We will also use $\w_i=(\w_{ij})_{j\in\Z^d}$.
On the space of environments $(\Omega,\mathfrak{S})$, we are given a 
certain $T$-invariant probability measure $\P$, with
$(\Omega,\mathfrak{S},(T^k)_{k\in\Z^d},\P)$ ergodic. 
We will say that the environment is i.i.d. when $\P$ is a product measure.
Let us now describe the process. First, the environment $\w$ is chosen
from the distribution $\P$. Once this is done, it remains fixed for
all times. The random walk in environment $\w$ is then the canonical
Markov chain $(X_n)_{n\geq0}$ with state space $\Z^d$ and
transition probability
\begin{eqnarray*}
&&P_x^\w(X_0=x)=1,\hfill\cr
&&P_x^\w(X_{n+1}=j\st X_n=i)=\pi_{ij}(\w).\hfill
\end{eqnarray*}
The process $P_x^\w$ is called the ``quenched law''. 
The ``annealed law'' is then
$$P_x=\int\!\!P_x^\w d\P(\w).$$
Already, one can see one of the difficulties of the model. When the
environment $\w$ is not fixed, i.e. under $P_0$, $X_n$ stops
being Markovian. 

Many questions arise about the different possible limit
theorems, such as the law of large numbers, central limit theorems, 
large deviation results, etc.
In the one-dimensional nearest-neighbor case,  the situation has been   
well understood, see e.g. \cite{sznitman,zeitouni} and the references 
therein. The reason for this is the possibility of explicit computations, and
the reversibility of the Markov chain.
In the higher-dimensional case, however, the amount of results is 
significantly less, see once again \cite{sznitman,zeitouni} for an overview. 

In the present paper, we are interested in the law of large numbers. In the
one-dimensional case, Solomon \cite{solomon} and Alili \cite{alili}
proved that the speed of escape of the particle (velocity at large
times) is a constant, $P_0$-a.s., that depends only on the distribution
of the environment. Later, Sznitman and Zerner~\cite{sz}, proved that,
under some technical transience condition on $\P$ (the so-called
Kalikow's condition), the law of large numbers still holds in the
multi-dimensional situation with i.i.d. environments.
To overcome the 
non-Markovian character of the walk, they used a renewal type
argument that appeared to be very specific to i.i.d. environments. 
Still, using the same method, Zeitouni \cite{zeitouni} proved
the law of large numbers, when i.i.d. environments are replaced by ones that 
are independent when a gap of size $L$ is allowed.
For more general mixing environments, the method seems to be too
rigid. However, physically relevant models, such as diffusions with random
coefficients, suggest that removing the independent environment
hypothesis is an important step towards a further understanding of
random walks in a random environment. For this, a different approach
is required. One approach that has proven to be very powerful in the
study of several other examples of random motions in random media, such
as in the works of Kipnis and Varadhan \cite{kv}, De Masi et al. \cite{demasi},
Olla \cite{olla}, and Papanicolaou and Varadhan \cite{pv}, is termed the
``point of view of the particle''. In this approach, one considers the
process $(T^{X_n}\w)$ of the environment as seen from the particle. 
This process is now a Markov process, with initial distribution
$\P$. The new inconvenience is that this Markov process has for its
state space 
the huge set $\Omega$. To apply the standard ergodic theorem, Kozlov
\cite{kozlov} showed that one needs to find an ergodic 
measure $\P_\infty$, that is invariant for the process $(T^{X_n}\w)$, 
and absolutely continuous relative to $\P$, see also lemma K 
below. This approach works perfectly
in the one-dimensional case, since one can compute 
$\P_\infty$ explicitly, see \cite{alili}. Moreover, in this case, one
does not need the i.i.d. hypothesis.  The hard problem though, is to 
find such a measure.  In the case of balanced walks, see
\cite{lawler}, one can prove the existence of such a measure, without 
actually computing it.
Even though in the two cases we mentioned above, the method
of the point of view of the particle did solve the problem,
it seems to have so far been of little help in the more general cases
of random walks in random environments.

As one will see in section \ref{sec_singular} below, one can not always
expect to be able to find an invariant measure $\P_\infty$ that is
absolutely continuous relative to $\P$, in all of $\Omega$. 
However, when studying walks that are 
transient in some direction $\ell\in\R^d-\{0\}$,
one expects the trajectories to stay in some half-plane
$H_k=\{x\in\Z^d:x.\ell\geq k\}$, for $k\leq0$.
In this paper, we further develop the approach of the point of view of the
particle, to be able to use it in the investigation of
higher-dimensional random walks in a non-necessarily i.i.d. random 
environment.     
In theorem \ref{invariant_multi}, we show that the conclusion of
Kozlov's lemma still holds if $\P_\infty$ is absolutely continuous
relative to $\P$, in every half-space $H_k$, instead of all of $\Z^d$. 
Then, in theorem \ref{multilln}, we show that Kalikow's condition implies
that, after having placed the walker at the origin,  the trajectories
do not spend ``too much'' time inside any half-plane $H_k$, collecting
therefore little information about the environment in there, and
satisfying the hypotheses of theorem \ref{invariant_multi}. This
implies the law of large numbers we are aiming for.

We will need the following definition.
We say that we have a finite range environment, or that the 
walk has finite range $M<\infty$, when
$$\P(\pi_{ij}=0~ ~ \mbox{when }\abs{i-j}>M)=1.$$
Throughout the rest of this work, we will only consider finite range
random walks in a random environment.

Let us now explain the structure of this paper. Section
\ref{sec_kalikow} introduces Kalikow's condition. There, we give an
effective condition that implies Kalikow's condition, even when $\P$
is not a product measure. By an effective condition, we mean a
condition that can be checked directly on the environment.

In section \ref{sec_1dlln}, we start with a warm up calculation. We consider
the one-dimensional finite range situation. We do not assume $\P$ to
be a product measure. In theorem \ref{1dlln}, we prove the law of large 
numbers in the one-dimensional non-nearest neighbor case, under
Kalikow's condition. 

In section \ref{sec_singular}, we explain why, in general, one can not use 
Kozlov's lemma in the multi-dimensional setup. 

In section \ref{sec_multiinv}, we prove theorem \ref{invariant_multi},
that extends Kozlov's lemma. We show that, in order to have a law of
large numbers, it is enough for the invariant measure $\P_\infty$ to
be absolutely continuous relative to $\P$ only in certain ``relevant''
parts of $\Z^d$. 

In section \ref{sec_dobrushin}, we introduce the Dobrushin-Shlosman
strong mixing condition.

In section \ref{sec_multilln}, we use theorem \ref{invariant_multi} to 
prove that Kalikow's condition implies the law of large numbers for
finite range random walks in a mixing random environment. This is our 
main theorem \ref{multilln}.

\mysection{Kalikow's condition.}
\label{sec_kalikow}
Let us start with a definition.  We define the drift $D$ to be
\begin{eqnarray}
\label{D}
D(\w)=E_0^\w(X_1)=\sum_i i\pi_{0i}(\w).
\end{eqnarray}
When studying the law of large numbers, one could try to examine 
first the case when the environment satisfies some condition that 
guarantees a strong drift in some direction
$\ell\in\R^d-\{0\}$.  
One such condition was introduced by Kalikow \cite{kalikow}. 
\begin{eqnarray}
\label{kalikow}
\inf_{U\in\mathscr U}\inf_{x\in U}{\ds
E_0\left(\sum_{j=0}^{T_U}\ind(X_j=x)D(T^{X_j}\w).\ell\right)\over
\ds E_0\left(\sum_{j=0}^{T_U}\ind(X_j=x)\right)}=\eps>0,
\end{eqnarray}
where $T_U=\inf\{j\geq0:\,X_j\not\in U\}$, and $\mathscr U$ is 
ranging over all finite sets that contain $0$, and have a path of 
range $M$ passing through all its points. We will call such sets
$M$-connected.
The expectations involved in
the above condition are all finite and positive 
(\cf \cite[p. 756--757]{kalikow}), if one assumes 
the following ellipticity condition to hold:
\begin{eqnarray}
\label{ellipticity}
\matrix{
\hbox{There exists }\kappa(\P)\in(0,1)~ \hbox{such that}\hfill\cr
\P(\pi_{i,j}>\kappa~ \hbox{when }\abs{i-j}\leq M)=1.}
\end{eqnarray}
In some situations, we will assume, instead, the weaker ellipticity condition
\begin{eqnarray}
\label{weakellip}
\P(\forall j\mbox{ s.t. }j.\ell\geq0\mbox{ and }\abs{j}=1:\pi_{0j}>0)=1.
\end{eqnarray}
In the rest of this work, we will consider condition
(\ref{ellipticity}) to be part of Kalikow's condition (\ref{kalikow}).
Sznitman and Zerner's \cite{sz} law of large numbers was established  
under condition (\ref{kalikow}). As a matter of fact, Kalikow's condition,
in the one-dimensional i.i.d. nearest-neighbor case, is equivalent 
to the condition $\E(\rho)<1$, \cf \cite[p. 1866--1867]{sz}. 
According to Solomon \cite{solomon}, this condition
characterizes the situation of walks with a positive speed of escape.
This is not the case in higher dimensions. In fact, Sznitman \cite{sT} proved that, 
in the i.i.d. case,
Kalikow's condition implies a strictly more general condition (the so-called $T^{'}$
 condition), which also implies a law of large numbers with a positive velocity.
One way to motivate Kalikow's condition is revealed by 
proposition 1 in \cite[p. 757--758]{kalikow}.

Of course (\ref{kalikow}) is not very practical, since it is not a
condition on the environment. 
Clearly, if one has a non-nestling environment, that is if there
exists a $\delta>0$ such that $\P(D.\ell\geq\delta)=1$, then 
(\ref{kalikow}) holds. In the nestling case, however, there is a
condition that is more concrete than (\ref{kalikow}), that implies 
it, and at the same time follows from many other interesting
conditions on the drift, such as: $\P(D.\ell<0)>0$, but there exists 
a constant $\delta>0$, such that
$\P(D.\ell\geq\delta)>C_\delta$ large enough.
It has already been established in
\cite[p. 759--760]{kalikow} and \cite[p. 36--37]{sznitman} that under the 
hypothesis that the environment is i.i.d., 
\begin{eqnarray}
\label{D+>D-}
\E(D.\ell^+)>\kappa^{-1}\E(D.\ell^-)
\end{eqnarray}
implies (\ref{kalikow}). 
In fact, one can relax the i.i.d. hypothesis
as follows. 
Let $\w_{\not x}=(\w_y)_{y\not=x}$, and define $Q_{\w_{\not x}}$ to be
the regular conditional probability, knowing $\w_{\not x}$, 
$Q_x$ be the marginal of $\w_x$, and $Q_{\not x}$ the marginal of
$\w_{\not x}$.
\begin{proposition}
Suppose that $Q_{\not 0}$-almost surely, $Q_{\w_{\not 0}}\ll Q_0$, and that
there exist two positive constants $A$ and $B$, such
that for $Q_0\otimes Q_{\not 0}$-almost every $\w=(\w_0,\w_{\not 0})$ one
has 
\begin{eqnarray}
\label{dependence}
0<A\leq h(\w_0,\w_{\not 0})={dQ_{\w_{\not 0}}\over dQ_0}(\w_0)\leq B<\infty.
\end{eqnarray}
Then the ellipticity condition (\ref{ellipticity}), along with 
\begin{eqnarray}
\label{D+>D-2}
\E(D.\ell^+)>\kappa^{-1}BA^{-1}\E(D.\ell^-),
\end{eqnarray}
imply Kalikow's condition (\ref{kalikow}).
\end{proposition}

\proof
Fix $U\subset\Z^d$, with $0\in U$. Define, for $\w\in\Omega$, $x,y\in\Z^d$,
\begin{eqnarray*}
f_{\w}(x)&=&P_0^{\w}(\exists k\in[0,T_U)\,:\,X_k=x),\cr
g_{\w}(x,y)&=&P_{x+y}^{\w}(X_k\not=x~ \forall k\in[0,T_U]).
\end{eqnarray*}
Note that we have, for $x\in U$,
$$P_x^{\w}(X_k\not=x~ \forall k\in(0,T_U])=
\sum_{\abs{y}\leq M}\pi_{x,x+y}(\w)g_{\w}(x,y).$$
Once $x$ is visited before exiting $U$, the number of returns to $x$,
up to time $T_U$, is geometrically distributed with the above failure
probability. Therefore, for $x\in U$, we have
$$E_0^{\w}\left(\sum_{j=0}^{T_U}\ind(X_j=x)\right)=
{f_{\w}(x)\over\ds\sum_{\abs{y}\leq M}g_{\w}(x,y)\pi_{x,x+y}(\w)}~,$$
where the numerator is exactly the probability of visiting $x$ at
least once. And, since $f_{\w}(x)$ and $g_{\w}(x,y)$
are $\sigma(\w_z;\,z\not=x)$-measurable, 
one has  
\begin{eqnarray*}
\qquad\quad\!\!\!\!&&\!\!\!\!\!\!\!\!\!\!\!\!\!\!\!\!\!\!\!\!\!\!
\!\!\!\!\!
\int\!\!{f_{\w}(x)D(T^x\w).\ell
\over\ds\sum_{\abs{y}\leq M}g_{\w}(x,y)\pi_{x,x+y}(\w)}d\P(\w)=
\int\!\!dQ_{\not x}(\w_{\not x})\int\!\!{f_{\w}(x)D(T^x\w).\ell\over
\ds\sum_{\abs{y}\leq M}g_{\w}(x,y)\pi_{x,x+y}(\w)}h(\w_x,\w_{\not x})
dQ_x(\w_x)\cr
&\geq&\!\!\!\!\int\!\!dQ_{\not x}(\w_{\not x})\int\!\!{f_{\w}(x)\over
\ds\max_{\abs{y}\leq M}g_{\w}(x,y)}(AD.\ell^+(T^x\w)-
\kappa^{-1}BD.\ell^-(T^x\w))dQ_x(\w_x)\cr
&=&\!\!\!\!\E(AD.\ell^+-\kappa^{-1}BD.\ell^-)
\int\!\!{f_{\w}(x)\over\ds\max_{\abs{y}\leq M}g_{\w}(x,y)}d\P(\w)\cr
&\geq&\!\!\!\!\kappa\E(AD.\ell^+-\kappa^{-1}BD.\ell^-)\int\!\!{f_{\w}(x)
\over\ds\sum_{\abs{y}\leq M}g_{\w}(x,y)\pi_{x,x+y}(\w)}d\P(\w).
\end{eqnarray*}
Which is Kalikow's condition, with 
$\eps=\kappa\E(AD.\ell^+-\kappa^{-1}BD.\ell^-)>0$.\qed\bigskip

Notice that in the i.i.d. case, (\ref{dependence}) clearly holds with 
$A=B=1$, and condition (\ref{D+>D-2}) is the same as (\ref{D+>D-}).
(\ref{dependence}) can also be easily checked, in the case of Gibbs
specifications, that we will use in the higher-dimensional case, see
section \ref{sec_dobrushin} below.
In this case, there exists a $C_1$
(same as in (\ref{radon_inf_bound})) such that the marginal $\mu_0$, of the
reference measure, satisfies
$$C_1^{-2}Q_0\leq C_1^{-1}\mu_0\leq Q_{\w_{\not 0}}\leq C_1\mu_0\leq 
C_1^2 Q_0.$$

Next, we show two implications of Kalikow's condition (\ref{kalikow}). 
Firstly, the walk has a ballistic character, in the following sense. 
\begin{lemma}
\label{kalikow_ballistic}
Assume we have a finite range environment for which
Kalikow's condition (\ref{kalikow}) holds. Let
$U\subset\Z^d$ be an $M$ connected set containing $0$, for which $E_0(T_U)<\infty$. Then
$E_0(X_{T_U}.\ell)\geq\eps E_0(T_U)$.
\end{lemma}

\proof For a finite $U$,  Kalikow's condition implies that
$$E_0\left(\sum_{j=0}^{T_U}\ind(X_j=x)D(T^{X_j}\w).\ell\right)\geq
\eps E_0\left(\sum_{j=0}^{T_U}\ind(X_j=x)\right).$$
Summing over all $x\in U$, and using that $D(T^{X_j}\w)=
E_0^\w(X_{j+1}-X_j|\mathscr{F}_j)$, and that $T_U$ is a stopping time,
one has
$$E_0\left(\sum_{j=0}^{T_U-1}(X_{j+1}-X_j).\ell\right)\geq
\eps E_0(T_U).$$
The claim follows. For an infinite $U$, the lemma follows from the
monotone convergence theorem, by taking 
increasing limits of finite sets.\qed\bigskip

The other consequence of Kalikow's condition is that,
under this condition, the walk almost surely escapes to
infinity in direction $\ell$. This was originally proved by Kalikow
\cite{kalikow}, and we reprove it here for the sake of completeness.
We also prove that the number of returns to the origin has a finite
annealed expectation. 
\begin{lemma}
\label{infinity}
Under Kalikow's condition (\ref{kalikow}), we have, 
\begin{eqnarray}
\label{infinity_eqn}
P_0\left(\lim_{n\rightarrow\infty}X_n.\ell=\infty\right)=1
\end{eqnarray}
and
\begin{eqnarray}
\label{visits}
\sum_{j\geq0}P_0(X_j=0)<\infty.
\end{eqnarray}
\end{lemma}

\proof
Let $U\subset\Z^d$ be a finite $M$-connected set containing $0$. 
Rewriting
(\ref{kalikow}), multiplying both sides by $e^{-\lambda x.\ell}$, 
for $\lambda>0$, and summing over all $x\in U$, one has  
\begin{eqnarray}
\label{lambdaa}
E_0\left(\sum_{j=0}^{T_U-1}e^{-\lambda X_j.\ell}D(T^{X_j}\w).\ell\right)
\geq\eps E_0\left(\sum_{j=0}^{T_U-1}e^{-\lambda X_j.\ell}\right).
\end{eqnarray}
On the other hand, since $T_U$ is a stopping time, one can write
$$\sum_{j=1}^{T_U}E_0^\w\left(e^{-\lambda X_j.\ell}\,|\,
\mathscr{F}_{j-1}\right)=
\sum_{j\geq1}E_0^\w\left(\ind(T_U\geq j)e^{-\lambda X_j.\ell}\,|\,
\mathscr{F}_{j-1}\right).$$
Hence, we have 
\begin{eqnarray*}
E_0\left(\sum_{j=1}^{T_U}e^{-\lambda X_j.\ell}\right)&=&
E_0\left(\sum_{j=1}^{T_U}E_0^\w\left(e^{-\lambda X_j.\ell}\,|\,
\mathscr{F}_{j-1}\right)\right)\cr
&=&E_0\left(\sum_{j=1}^{T_U}e^{-\lambda X_{j-1}.\ell}\left(1-\lambda
D(T^{X_{j-1}}\w).\ell+O(M^2\lambda^2)\right)\right)\cr
&\leq&E_0\left(\sum_{j=0}^{T_U-1}e^{-\lambda X_j.\ell}\right)
(1-\lambda\eps+O(M^2\lambda^2)).
\end{eqnarray*}
where we have used (\ref{lambdaa}) to get the inequality. 
Taking $\lambda>0$ small enough, and increasing $U$ to all of $\Z^d$,
one has
\begin{eqnarray}
\label{TU}
E_0\left(\sum_{j\geq0}e^{-\lambda X_j.\ell}\right)<\infty
\end{eqnarray}
and, therefore, 
$$P_0\left(\linf_{n\rightarrow\infty}X_n.\ell<\infty\right)\leq
P_0\left(\sum_{j\geq0}e^{-\lambda X_j.\ell}=\infty\right)=0,$$
proving (\ref{infinity_eqn}). 
Using (\ref{TU}), one also proves (\ref{visits})\bigskip

\hfill
$\displaystyle{\sum_{j\geq0}P_0(X_j=0)=E_0\left(\sum_{j\geq0}\ind(X_j=0)\right)
\leq E_0\left(\sum_{j\geq0}e^{-\lambda X_j}\right)<\infty.}$
\hfill\qed\bigskip

Next, as a warm up for the method we will use later in the
multi-dimensional situation, we examine the simpler case of
one-dimensional random walks. 

\mysection{The One-dimensional Case.}
\label{sec_1dlln}
In this section, we will prove the law of large numbers for
one-dimensional finite range random walks in a random environment. 
Let us recall a lemma, also valid for $d\geq2$, that
was proved by Kozlov in \cite{kozlov}. \bigskip

{\sc Lemma K. (Kozlov \cite{kozlov})} {\it
Assume that the weak ellipticity condition (\ref{weakellip}) holds. 
Suppose also that there exists a probability measure
$\P_\infty$ that is invariant for the process $(T^{X_n}\w)_{n\geq0}$,
and that is absolutely continuous relative to the ergodic $T$-invariant 
environment $\P$. Then
\begin{enumerate}
\item The measures $\P$ and $\P_\infty$ are in fact mutually absolutely continuous.
\item The Markov process $(T^{X_n}\w)_{n\geq0}$ with initial distribution 
$\P_\infty$ is ergodic.
\item There can be at most one such $\P_\infty$.  
\item The following law of large numbers is satisfied.
$$P_0\left(\lim_{n\rightarrow\infty}{X_n\over n}=\E^{\P_\infty}(D)\right)=1,$$
\end{enumerate}
where $D$ is the drift defined in (\ref{D}).
}\bigskip

One then has the following theorem.
\begin{theorem}
\label{1dlln}
Under Kalikow's condition (\ref{kalikow}), with $\ell=1$, the process
$(T^{X_n}\w)_{n\geq0}$ has an invariant measure $\P_\infty$ that is
absolutely continuous relative to $\P$, and we have a law of large
numbers for finite range random walks in the ergodic $T$-invariant
environment $\P$. 
$$\P\left(\lim_{n\rightarrow\infty}{X_n\over n}=
\E^{\P_\infty}(D)\right)=1.$$
\end{theorem}

\proof
Define
$$g_{ij}(\w)=\sum_{n\geq0}P_i^\w(X_n=j)=E_i^\w(N_j),$$
where $N_j$ is the number of visits of the random walk to site $j$. 
The renewal property gives, for $i\not=j$,
$$g_{ij}=E^\w_j(N_j)P^\w_i(V_j<\infty)\leq g_{jj},$$
with $V_j=\inf\{n>0\,:\,X_n=j\}$. Moreover,  $g_{jj}$'s are all
identically distributed in the annealed setting. Thus, according to
(\ref{visits}), they are all in $L^1(\Omega,\P)$. 
For $i\leq j$ define 
$$G_{ij}={1\over j-i+1}\sum_{k=i}^j g_{kj}\leq g_{jj}.$$
Using the diagonal trick, one can extract a subsequence of the 
$G_{ij}$'s that converges weakly, as $i$ decays to $-\infty$, to a limit 
$\mu_j\in L^1(\Omega,\P)$, for all $j$. Then, for any fixed $j$,
$\mu_j$ is a limit point for the $g_{ij}$'s as well. Using the diagonal 
trick again, one can find a subsequence of the $g_{ij}$'s, that converges 
weakly to $\mu_j$, for all $j$.
We will keep referring to both subsequences by $G_{ij}$ and $g_{ij}$.

Notice that if $k\not=j$, then
$$\sum_i\pi_{ij}g_{ki}=\sum_{n\geq0}\sum_i\pi_{ij}P_k^\w(X_n=i)
=\sum_{n\geq0}P_k^\w(X_{n+1}=j)=g_{kj}.$$
Therefore, for $\P$-a.e. $\w$,
$\sum_i\pi_{ij}\mu_i=\mu_j$. Also, 
$$g_{i0}\circ T=\sum_{n\geq0}P^{T\w}_i(X_n=0)=
\sum_{n\geq0}P^\w_{i+1}(X_n=1)=g_{i+1,1},$$
and, the same holds for the $G_{ij}$'s. Therefore,
for $\P$-a.e. $\w$,
$$\mu_0(T\w)=\lim_{i\rightarrow-\infty}G_{i0}(T\w)
=\lim_{i\rightarrow-\infty}G_{i+1,1}(\w)=\mu_1(\w).$$
This shows that $\mu_0d\P$ is an invariant measure
for the process $(T^{X_n}\w)_{n\geq0}$. 
Next, we need to show that $\mu_0$ is not trivial. To this end, we
recall lemma \ref{infinity}. According to this lemma, Kalikow's
condition implies that, for $\P$-a.e. $\w$,
$P_0^\w\left(\lim_{n\rightarrow\infty}X_n=\infty\right)=1$.
The finite range character of the walk implies then that for each
$i<j$, $\P$-a.s., $\sum_{k=j}^{j+M-1}g_{ik}\geq1$.
Taking 
the limit in $i$, we have that, $\P$-a.s.,
$\sum_{k=j}^{j+M-1}\mu_k\geq1$. Therefore, by the ergodicity of $\P$,
$\E(\mu_0)\geq M^{-1}$.

Defining $\P_\infty$ such that 
$${d\P_\infty\over d\P}={\mu_0\over\E(\mu_0)}$$ 
gives an invariant probability measure for the process of the
environment, as seen from the particle. This measure is absolutely
continuous relative to $\P$, and lemma K 
concludes the proof.\qed\bigskip

Now, we move to the multi-dimensional situation. In the following
section, we will show why it is quite different from the situation
above, and why Kozlov's lemma K 
can not be used. 

\mysection{Motivation.}
\label{sec_singular}
Consider the case where $d=2$, the environment is i.i.d.,
and $$\P(\pi_{(0,0)(1,0)}=1)=\P(\pi_{(0,0)(0,1)}=1)=0.5.$$
Once the environment is chosen, the walk is determined,
following the assigned directions.
The annealed process is in fact the same as $0.5(n-S_n,n+S_n)$, with 
$S_n$ a one-dimensional simple symmetric random walk.
Therefore, one obviously has the following law of large numbers:
$$P_0\left(\lim_{n\rightarrow\infty}{X_n\over n}=(0.5,0.5)\right)=1.$$
Yet, defining $\P_n$, to be the measure on the environment as seen
from the particle at time $n$:
$$\P_n(A)=P_0(T^{X_n}\w\in A)$$
and ${\mathfrak S}_{-k}$, as the $\sigma$-algebra
generated by the environment at sites $x$ such that
$x.(1,1)\geq-k$, one has the following proposition.
\begin{proposition}
\label{singular} There exists a probability measure $\P_\infty$ to 
which $\P_n$ converges weakly. Moreover, $\P_\infty$ is mutually singular 
with $\P$, and there is no probability measure that is at the 
same time,
invariant for $(T^{X_n}\w)$, and absolutely continuous relative to $\P$.
Furthermore, for $k\leq n$, one has ${\P_n}_{\big|{\mathfrak S}_{-k}}=
{\P_k}_{\big|{\mathfrak S}_{-k}}$, and therefore,
${\P_\infty}_{\big|{\mathfrak S}_{-k}}=
{\P_k}_{\big|{\mathfrak S}_{-k}}\ll\P_{\big|{\mathfrak S}_{-k}}$.
\end{proposition}

For a complete proof, see propositions 1.4. and 1.5. in \cite{dynstat}.
Although the ellipticity condition is not satisfied, this
model is instructive.  It shows us that to prove a law of large
numbers, one need not necessarily look for a $\P_\infty$ that is
absolutely continuous relative to $\P$ on the whole space. Instead,
maybe one should try to prove that $\P_\infty\ll\P$ in the
``relevant'' part of the space, that is all half-spaces 
$\{x:x.(1,1)\geq-k\}$, for $k\geq0$.
This is still much weaker than absolute continuity in
the whole space. We will address this issue, in the following section.

\mysection{On the invariant measure for {\boldmath $d\geq2$}.}
\label{sec_multiinv}
For $k\in\Z$, let ${\mathfrak S}_k=\sigma(\w_x\,:\,x.\ell\geq k)$ be 
the $\sigma$-algebra generated by the part of the environment in the
right half-plane $H_k=\{x\,:\,x.\ell\geq k\}$. In this section, 
we will not assume the ellipticity condition (\ref{ellipticity})
to hold. Instead, we will assume the weaker ellipticity condition 
(\ref{weakellip}) we assumed in lemma K. 
We modify lemma K, as suggested by the example in section \ref{sec_singular},
and we have the following theorem.
\begin{theorem}
\label{invariant_multi}
Let $\P$ be ergodic, and $T$-invariant, with finite range $M$.
Assume that the weak ellipticity condition (\ref{weakellip}) holds, and that 
\begin{eqnarray}
\label{assumption}
P_0\left(\lim_{n\rightarrow\infty}X_n.\ell=\infty\right)=1.
\end{eqnarray}
Suppose also that there exists a probability measure
$\P_\infty$ that is invariant for the process $(T^{X_n}\w)_{n\geq0}$,
and that is absolutely continuous relative to $\P$, in every
half-space $H_k$, with  $k\leq0$. Then
\begin{enumerate}
\item The measures $\P$ and $\P_\infty$ are in fact mutually absolutely
continuous on every $H_k$, with $k\leq0$. 
\item The Markov process $(T^{X_n}\w)_{n\geq0}$ with initial distribution 
$\P_\infty$ is ergodic.
\item There can be at most one such $\P_\infty$, and if
$\tilde\P_n(A)=n^{-1}\sum_{m=1}^n P_0(T^{X_m}\w\in A)$, then
$\tilde\P_n$ converges weakly to $\P_\infty$.
\item The following law of large numbers is satisfied.
$$\P\left(\lim_{n\rightarrow\infty}{X_n\over n}=\E^{\P_\infty}(D)\right)=1.$$
\end{enumerate}
\end{theorem}

\proof\medskip

$\bullet$ $\left(\forall k\leq0:\,{\P_\infty}_{\big|{\mathfrak S}_k}
\sim\P_{\big|{\mathfrak S}_k}\right)$:
Fix $k\leq0$, and let 
$G_k={{d\P_\infty}_{|{\mathfrak S}_k}\over
d\P_{|{\mathfrak S}_k}}$.
Then
\begin{eqnarray*}
0&=&\int_{\{G_k=0\}}\!\!\!\!G_k d\P=\int\!\!\ind_{\{G_k=0\}} d\P_\infty
=\int\!\!\sum_{\abs{e}\leq M}\pi_{0e}\,\ind_{\{G_k=0\}}\circ T^e\,
d\P_\infty\cr
&\geq&\int\!\!
\sum_{{}^{\abs{e}=1}_{e.\ell\geq0}}\pi_{0e}\,
\ind_{\{G_k=0\}}\circ T^e\,G_k\,d\P
=\int_{\{G_k=0\}}\!\sum_{{}^{\abs{e}=1}_{e.\ell\leq0}}
\pi_{e0}G_k\circ T^{e}\,d\P,
\end{eqnarray*}
where the inequality used the fact that if $e.\ell\geq0$ then
$G_k\circ T^e$ is still ${\mathfrak S}_k$-measurable. 
Using the weak ellipticity condition (\ref{weakellip}), the above 
inequality implies that $\P$-a.s. we have
$\{G_k=0\}\subset T^e\{G_k=0\}$,
when $\abs{e}=1$ and $e.\ell\geq0$. Since $T$ is $\P$-preserving, we
have
$\{G_k=0\}=T^e\{G_k=0\}$, $\P$-a.s.
And since $(T^e)_{{}^{\abs{e}=1}_{e.\ell\geq0}}$ generates the group
$(T^x)_{x\in\Z^d}$, we have that $\{G_k=0\}$ is
$\P$-a.s. shift-invariant. But $\P$ is ergodic, and thus 
$\P(G_k=0)$ is $0$ or $1$. However, $\E(G_k)=1$, and therefore 
$\P(G_k>0)=1$, and
$\P_\infty$ and $\P$ are mutually absolutely continuous on $H_k$, for
any $k\leq0$.\medskip

$\bullet$ Ergodicity of $(T^{X_n}\w)_{n\geq0}$ with initial
distribution $\P_\infty$:
Consider a bounded local function $f$ on $\Omega$ that is ${\mathfrak
S}_K$-measurable, for some $K\leq0$. Define 
$g=\E^{\P_\infty}(f|{\mathcal I})$, where $\mathcal I$
is the invariant $\sigma$-field for the process $(T^{X_n}\w)_{n\geq0}$.
Birkhoff's ergodic theorem implies that for $\P_\infty$-a.e. $\w$
\begin{eqnarray}
\label{birkhoff}
P^\w_0\left(\lim_{n\rightarrow\infty}
n^{-1}\sum_{m=1}^n f(T^{X_m}\w)=g(\w)\right)=1.
\end{eqnarray}
Using the fact that $\P_\infty$ is invariant and that $g$ is harmonic,
we have
\begin{eqnarray*}
\sum_{\abs{e}\leq M}\int\!\!\pi_{0e}(g-g\circ T^e)^2\,d\P_\infty&=&
\int\!\!g^2\,d\P_\infty-2\int\!\!g\sum_{\abs{e}\leq M}\pi_{0e}\,
g\circ T^e\,d\P_\infty\cr
&&\quad+\int\!\!\sum_{\abs{e}\leq M}\pi_{0e}(g\circ T^e)^2\,d\P_\infty=0.
\end{eqnarray*}
Noticing that $\pi_{0e}$ is ${\mathfrak S}_0$-measurable we conclude
that the above equation, along with the weak ellipticity condition
(\ref{weakellip}), implies that for $\abs{e}=1$ and $e.\ell\geq0$
\begin{eqnarray}
\label{g=gt}
g=g\circ T^e\qquad \P_\infty\mbox{-a.s.} 
\end{eqnarray}

Moreover, if we define 
$$S=\left\{\w:\,\forall y\in\Z^d~ ~ 
P_y^\w\left(\inf_{m\geq0}X_m.\ell<0\right)=1\right\},$$
then $\P(S)=0$. This is because otherwise the renewal property for the
quenched walk would imply that $P_0(X_n.\ell<0\mbox{ i.o.})>0$, and this
contradicts (\ref{assumption}). Hence, we have that for $\P$-a.e. $\w$,
there exists a $y$ such that
$$P_y^\w\left(\inf_{m\geq0}X_m.\ell\geq0\right)>0.$$
In particular, $y.\ell\geq0$.
The weak ellipticity condition (\ref{weakellip}) implies that
for $\P$-a.e. choice of $\w$, the walk starting at $0$ will,
with positive probability under $P_0^\w$, reach $y$ without
backtracking below $0$. This means that
$$P_0^\w\left(\inf_{m\geq0}X_m.\ell\geq0\right)>0\qquad \P\mbox{-a.s.}$$
But the above event is ${\mathfrak S}_0$-measurable, and therefore we
have
$$P_0^\w\left(\inf_{m\geq0}X_m.\ell\geq0\right)>0\qquad 
\P_\infty\mbox{-a.s.}$$

Define now
$${\bar g}(\w)=P_0^\w\left(\inf_{m\geq0}X_m.\ell\geq0\right)^{-1}
\lsup_{n\rightarrow\infty}\int_{\{\inf_{m\geq0}X_m.\ell\geq0\}}
n^{-1}\sum_{m=1}^n f(T^{X_m}\w)dP_0^\w.$$
Then because of (\ref{birkhoff}), we know that
$g={\bar g},~ \P_\infty$-a.s. However,
it is clear that $\bar g$ is ${\mathfrak S}_K$-measurable.
Formula (\ref{g=gt}) then implies that $g=g\circ T^e,~ \P$-a.s. and the
ergodicity of $\P$ implies that $g$ is constant $\P$-a.s., and thus 
$\P_\infty$-a.s. This proves that the invariant $\sigma$-field
$\mathcal I$ is trivial, and that concludes the proof of ergodicity of 
$(T^{X_n}\w)_{n\geq0}$ with initial distribution $\P_\infty$.\medskip

$\bullet$ Uniqueness of $\P_\infty$:
Let $f$ be a local bounded ${\mathfrak S}_K$-measurable function, 
for $K\leq0$. Notice that due to ergodicity, we have $\P_\infty$-a.s.
$$E_0^\w\left(\lim_{n\rightarrow\infty}
n^{-1}\sum_{m=1}^n f(T^{X_m}\w)\right)=\E^{\P_\infty}(f)$$
and, therefore, for $k\leq0$, we have $\P_\infty$-a.s.
$$E_0^\w\left(\lim_{n\rightarrow\infty}
n^{-1}\sum_{m=1}^n f(T^{X_m}\w);
\inf_{m\geq0}X_m.\ell\geq k\right)=\E^{\P_\infty}(f)
P_0^\w\left(\inf_{m\geq0}X_m.\ell\geq k\right).$$
Both functions above are ${\mathfrak S}_{k+K}$-measurable. 
Therefore, the same equation holds $\P$-a.s.
Integrating over $\w$,  one has
\begin{eqnarray*}
\E^{\P_\infty}(f)&=&\lim_{k\rightarrow-\infty}
E_0\left(\lim_{n\rightarrow\infty}
n^{-1}\sum_{m=1}^n f(T^{X_m}\w);
\inf_{m\geq0}X_m.\ell\geq k\right)\cr
&=&\lim_{n\rightarrow\infty}n^{-1}\sum_{m=1}^n E_0(f(T^{X_m}\w))
=\lim_{n\rightarrow\infty}\E^{\tilde\P_n}(f),
\end{eqnarray*}
which uniquely defines $\P_\infty$ as the weak limit of $\tilde\P_n$.\medskip

$\bullet$ The law of large numbers:
Taking $f$ to be the drift $D$, we have, for all $k\leq0$, and 
$\P_\infty$-a.e. $\w$
$$P_0^\w\left(\lim_{n\rightarrow\infty}
n^{-1}\sum_{m=1}^n D(T^{X_m}\w)=\E^{\P_\infty}(D);
\inf_{m\geq0}X_m.\ell\geq k\right)=
P_0^\w\left(\inf_{m\geq0}X_m.\ell\geq k\right).$$
Once again, this is also true $\P$-a.s. and taking $k$ to $-\infty$ 
we have
\begin{eqnarray}
\label{lln_drift}
P_0\left(\lim_{n\rightarrow\infty}
n^{-1}\sum_{m=1}^n D(T^{X_m}\w)=\E^{\P_\infty}(D)\right)=1.
\end{eqnarray}
For the rest of the proof, we follow the argument in
\cite[p.~10]{sznitman}. To this end, 
$M_n=X_n-X_0-\sum_{m=0}^{n-1} D(T^{X_m}\w)$ is a martingale
with bounded increments under $P_0^\w$. Therefore
$P_0^\w\left(\lim_{n\rightarrow\infty}{n^{-1}M_n}=0\right)=1$.
Combining this with (\ref{lln_drift}), one obtains the desired law of
large numbers.\qed\bigskip

Next, we will relax the absolute continuity condition to a weaker, but
sufficient, condition.
But first, we need some definitions. For a measure $\P_\infty$, and $k\leq0$, 
define $\P_\infty^{k,\ll}$ (resp. $\P_\infty^{k,\perp}$) 
to be the absolutely continuous (resp. singular) part of 
${\P_\infty}_{\big|{\mathfrak{S}_k}}$ relative to $\P_{\big|{\mathfrak{S}_k}}$. 
For $A\in{\mathfrak{S}_k}$, and $j\leq k$, $\P_\infty^{j,\ll}(A)$ (resp. 
$\P_\infty^{j,\perp}(A)$) is a monotone sequence, and
there exists a measure $\P_\infty^{\infty,\ll}$ (resp. 
$\P_\infty^{\infty,\perp}$) such that, 
$\P_\infty^{\infty,\ll}(A)=\inf_{j\leq k}\P_\infty^{j,\ll}(A)$ (resp.
$\P_\infty^{\infty,\perp}(A)=\sup_{j\leq k}\P_\infty^{j,\perp}(A)
=\P_\infty(A)-\P_\infty^{\infty,\ll}(A)$).
Now, we have the following lemma.
\begin{lemma}
\label{invariant_abs2} If $\P_\infty$ is invariant for the process
$(T^{X_n}\w)_{n\geq0}$ and $\P_\infty^{\infty,\ll}(\Omega)>0$, then
$\hat\P_\infty=\P_\infty^{\infty,\ll}(\Omega)^{-1}\P_\infty^{\infty,\ll}$ is a 
probability
measure that is also invariant. Moreover, $\hat\P_\infty$ is absolutely
continuous relative to $\P$, in every half-space $H_k$, with  $k\leq0$.
\end{lemma}

\proof 
One clearly has ${\P_\infty^{\infty,\ll}}_{\big|{\mathfrak{S}_k}}\leq
\P_\infty^{k,\ll}\ll\P_{\big|{\mathfrak{S}_k}}$. This proves the absolute continuity 
part of the claim of the lemma. To show the invariance of $\hat\P_\infty$, it is
enough to show the invariance of $\P_\infty^{\infty,\ll}$. To this end,
denote the transition probability of the process of the 
environment viewed from the particle by 
$$\pi(\w,A)=\sum_{\abs{e}\leq M}\pi_{0e}(\w)\ind_A(T^e\w),$$
and define the operator $\Pi$, acting on measures, as
$$\Pi\P(A)=\int\!\!\pi(\w,A)d\P(\w).$$
Now, consider $A\in\mathfrak{S}_{k+M}$, with $\P(A)=0$. Since $\Pi\P\ll\P$, we
have $\Pi\P(A)=0$. Therefore, $\pi(\w,A)=0$, 
$\P_{\big|{\mathfrak{S}_k}}$-a.s. and thus $\P_\infty^{k,\ll}$-a.s.
as well. Hence, $\Pi\P_\infty^{k,\ll}(A)=0$. This proves that
$\Pi\P_\infty^{k,\ll}\ll\P_{\big|{\mathfrak{S}_{k+M}}}$
and, since $\Pi\P_\infty^{k,\ll}\leq\Pi({\P_\infty}_{\big|{\mathfrak{S}_k}})=
{\P_\infty}_{\big|{\mathfrak{S}_{k+M}}}$, we have
$\Pi\P_\infty^{k,\ll}\leq\P_\infty^{k+M,\ll}$. 
Taking limits,
one has $$\Pi\P_\infty^{\infty,\ll}\leq\P_\infty^{\infty,\ll}.$$
However, the two measures above give the same mass to $\Omega$, and therefore are
equal.\qed
\begin{remark}
Given an invariant measure $\P_\infty$, one can decompose it, relative to
$\P$, into $\P_\infty^\ll$ and $\P_\infty^\perp$. Using the same argument as 
above, it is easy to see that $\P_\infty^\ll$
is again invariant, and that $\P_\infty^\ll\leq\P_\infty^{\infty,\ll}$. Due 
to the uniqueness of the measure in theorem \ref{invariant_multi}, one sees 
that if $\P_\infty^\ll$ is not trivial, then $\P_\infty^\ll$ and 
$\P_\infty^{\infty,\ll}$ are proportional. Therefore, the latter is absolutely
continuous, relative to $\P$, in the whole space, and thus 
$\P_\infty^{\infty,\ll}\leq\P_\infty^\ll$, and
$\P_\infty^\ll=\P_\infty^{\infty,\ll}$.
\end{remark}
Before we move to the discussion of the law of large numbers, we will
introduce, and recall some facts about the Dobrushin-Shlosman mixing
condition for random fields. 

\mysection{The Dobrushin-Shlosman mixing condition.}
\label{sec_dobrushin}
First, we introduce some notations.
For a set $V\subset\Z^d$, let us denote by $\Omega_V$ the set of
possible configurations $\w_V=(\w_x)_{x\in V}$, and by $\mathfrak{S}_V$ the
$\sigma$-field generated by the environments $(\w_x)_{x\in V}$. 
For a probability measure $\P$, we will denote by $\P_V$, 
the projection of $\P$ onto $(\Omega_V,\mathfrak{S}_V)$. 
For $\w\in\Omega$, denote by 
$\P_V^\w$ the regular conditional probability, knowing
$\mathfrak{S}_{\Z^d-V}$, on $(\Omega_V,\mathfrak{S}_V)$.
Furthermore, for $\Lambda\subset V$, $\P_{V,\Lambda}^\w$ will denote
the projection of $\P_V^\w$ onto $(\Omega_\Lambda,\mathfrak{S}_\Lambda)$.
Also, we will use the notation $V^c=\Z^d-V$,
$\partial_r V=\{x\in\Z^d-V:\,\mbox{dist}(x,V)\leq r\}$, with $r\geq0$,
and $\card{V}$ will denote the cardinality of $V$.
Finally, for $\w,\bar\w\in\Omega$, $V,W\subset\Z^d$ with 
$V\cap W=\emptyset$, we will use $(\bar\w_V,\w_W)$ to denote
$\bar{\bar\w}_{V\cup W}$ such that $\bar{\bar\w}_V=\bar\w_V$ and 
$\bar{\bar\w}_W=\w_W$.
We will also need the following definitions. 

By an $r$-specification ($r\geq0$) we mean a system of 
functions $Q=\{Q_V^\cdot(\cdot):V\subset\Z^d,\card{V}<\infty\}$, 
such that for all $\w\in\Omega$, $Q_V^\w$ is a probability measure on 
$(\Omega_V,\mathfrak{S}_V)$, and, for all $A\in\mathfrak{S}_V$, 
$Q_V^\cdot(A)$ is $\mathfrak{S}_{\partial_rV}$-measurable.
Sometimes, for notational convenience, $Q_V^\cdot(A)$ will be thought 
of as a function on $\Omega_{\partial_rV}$. For $\Lambda\subset V$,
we will denote by $Q_{V,\Lambda}^\w$, the projection of $Q_V^\w$ onto 
$(\Omega_\Lambda,\mathfrak{S}_\Lambda)$.

A specification $Q$ is self-consistent if, for any finite $\Lambda,V$,
$\Lambda\subset V\subset\Z^d$, one has, for $Q_V^\w$-a.e. $\bar\w_V$,
$(Q_V^\w)_\Lambda^{\bar\w_V}=Q^{(\w_{V^c},\bar\w_V)}_\Lambda$.
We will say that a probability measure $\P$ is consistent
with a specification $Q$, if $\P_V^\w$ coincides with 
$Q_V^\w$, for every finite $V\subset\Z^d$ and $\P$-a.e. $\w$.
Notice that this can only happen when $Q$ is self-consistent. 
In this case, $Q$ is uniquely determined by $\P$. 
The question is, however, whether $Q$ determines $\P$, and whether it
does so uniquely. To this end, Dobrushin and Shlosman \cite{ds} gave
a sufficient condition to answer the above questions positively.\bigskip

{\sc Theorem DS. (Dobrushin-Shlosman \cite{ds})} 
{\it Let $Q$ be a self-consistent
$r$-specification, and assume the Dobrushin-Shlosman strong decay 
property holds, i.e. there exist $G,g>0$ such that for all 
$\Lambda\subset V\subset\Z^d$ finite, $x\in\partial_r V$, 
and $\w,\bar\w\in\Omega$, 
such that $\w_y=\bar\w_y$ when $y\not=x$, we have
\begin{eqnarray}
\label{mixing}
\mbox{Var}(Q_{V,\Lambda}^\omega,Q_{V,\Lambda}^{\bar\omega})\leq 
Ge^{-g\,{\scriptstyle\mathrm{dist}}(x,\Lambda)},
\end{eqnarray}
where $\mbox{Var}(\cdot,\cdot)$ is the variational distance
$\mbox{Var}(\mu,\nu)=\ds\sup_{E\in\mathfrak{S}}(\mu(E)-\nu(E))$.
Then, there exists a unique $\P$ that is consistent with $Q$.
Moreover, we have, for all $\w\in\Omega$,
\begin{eqnarray}
\label{phase}
\lim_{\mathrm{dist}(\Lambda,V^c)\rightarrow\infty}
\mbox{Var}(Q_{V,\Lambda}^\omega,\P_\Lambda)=0.
\end{eqnarray}
}
\indent The main example of self-consistent specifications are Gibbs 
specifications. For the precise definition of a Gibbs specification
with inverse temperature $\beta>0$, see \cite{ds}. Moreover, if the 
interaction is translation-invariant, and the specification satisfies
(\ref{mixing}), then the unique field $\P$ is also shift-invariant;
see \cite[sec. 5.2]{hog}.
One should note that the conditions of theorem DS are satisfied when one
considers Gibbs fields in the high-temperature region, i.e. when $\beta$ 
is small; see \cite{ds_analiticity}.

We will need the following lemma. The proof depends on another
lemma, and will be outlined in the appendix at the end of the paper.
\begin{lemma}
\label{dobrushin_bound}
Let $(\P^\w_V)$ be a Gibbs $r$-specification satisfying
(\ref{mixing}), and let $\P$ be the unique translation-invariant Gibbs
field, consistent with $(\P^\w_V)$. 
Consider $H\subset\Z^d$, and $\Lambda\subset H^c$ with 
$\mbox{\rm dist}(\Lambda,H)>r$. Then
$$\sup_{F\in\mathfrak{F}}\sup_\w
{\E(F|\mathfrak{S}_H)(\w)\over\E(F)}\leq
\exp\left(C\sum_{{}^{x\in\partial_r(H^c)}_{y\in\partial_r(\Lambda^c)}}
e^{-g\,{\scriptstyle\mathrm{dist}}(x,y)}\right),$$
where
$\mathfrak{F}=\{F\geq0,\mathfrak{S}_{\Lambda}\mbox{-measurable},
\hbox{ s.t. }\E(F)>0\}$.
\end{lemma}

\mysection{The law of large numbers.}
\label{sec_multilln}
We need now to find an invariant measure $\P_\infty$ that is
absolutely continuous relative to $\P$, in each half-plane. The reason
for which such a measure would exist is a strong enough transience condition.
We will consider an environment that either satisfies the
Dobrushin-Shlosman mixing condition (\ref{mixing}), or that is
$L$-dependent in direction $\ell$, i.e. there exists $L>0$, such that
\begin{eqnarray}
\label{gap}
\sigma(\w_x;\,x.\ell\leq0)~ \hbox{and }\sigma(\w_x;\,x.\ell\geq L)
~ \hbox{are independent.}
\end{eqnarray}
The following is our main theorem.
\begin{theorem}
\label{multilln}
Suppose that $\P$ is of finite range, $T$-invariant, ergodic, and 
satisfies one of the mixing conditions (\ref{mixing}) or (\ref{gap}). 
Suppose also that the strong $\kappa$-ellipticity 
condition (\ref{ellipticity}) holds, and that Kalikow's condition 
(\ref{kalikow}), in direction $\ell\in S^{d-1}$, is satisfied.
Then, the process $(T^{X_n}\w)_{n\geq0}$ admits an invariant
probability measure  $\hat\P_\infty$ that is
absolutely continuous relative to $\P$, in every half-space $H_k$ 
with $k\leq0$, 
and we have a law of large numbers for the finite range 
random walk in environment $\P$, with a non-zero limiting velocity:
$$\P\left(\lim_{n\rightarrow\infty}{X_n\over n}=
\E^{\hat\P_\infty}(D)\not=0\right)=1.$$
Moreover, if $\P_n(A)=P_0(T^{X_n}\w\in A)$, i.e. $\P_n$ is the measure on the
environment as seen from the particle at time $n$, then 
$N^{-1}\sum_{n=1}^N\P_n$ converges weakly to $\hat\P_\infty$.
\end{theorem}

\proof Define the spaces 
$$\VV_n=\{\hbox{paths }\vv\hbox{, of range }M\hbox{, length }
n+1\hbox{, and ending at }0\}$$
and the space $\VV$ of paths $\vv$, of range $M$, ending at
$0$, and of either finite or infinite length. 
Being a closed subspace of 
$(\{e\in\Z^d:\abs{e}\leq M\}\cup\{\hbox{`Stop'}\})^\N$, 
endowed with the product topology, $\VV$ 
is compact. And, if we now consider the space $\VV_\infty\subset\VV$ of 
paths of range $M$ and of infinite length, that end at $0$, then $\VV_\infty$ 
is again a compact space.

Let us now define a sequence of measures $R_n$ on
$\VV\times\Omega$ as follows. Clearly, $R_n$ will be supported on
$W_n\times\Omega$, and for $\vv=(x_0,x_1,x_2,\cdots,x_n=0)\in W_n$,
$A\in\mathfrak{S}$, 
$$R_n(\{\vv\}\times A)=
P_0((-X_n,X_1-X_n,\cdots,X_{n-1}-X_n,0)=\vv,T^{X_n}\w\in A).$$
Notice that $\P_n$ is the marginal of $R_n$, and therefore, the
disintegration lemma implies that
$$\P_n(A)=\int\!\!\P_\svv(A)dQ_n(\vv),$$
where $Q_n$ is the marginal of $R_n$ over $W_n$. It assigns 
probability $\E(\pi_\svv)$ to paths $\vv$ of length $n+1$, and 
ending at $0$. Here, $$\pi_\svv=\prod_{i=0}^{n-1}\pi_{x_i,x_{i+1}}.$$
In fact, one can compute $\P_\svv$ explicitly. Indeed, 
\begin{eqnarray*}
\P_n(A)&=&P_0(T^{X_n}\w\in A)=\sum_{x\in\Z^d}P_0(X_n=x,T^x\w\in A)
=\sum_{x\in\Z^d}P_x(X_n=0,\w\in A)\cr
&=&\int_A\!\sum_{x\in\Z^d}P_x^\w(X_n=0)d\P(\w)
=\int_A\sum_{\svv\in\sVV_n}\pi_\svv(\w)d\P(\w).
\end{eqnarray*}
Using Fubini's theorem, we have
$$\P_n(A)=\int\!\!\P_\svv(A)dQ_n(\vv),\quad\hbox{with}\qquad
{d\P_\svv\over d\P}={\pi_\svv\over\E(\pi_\svv)}.$$
The measure $\P_\svv$ could be thought of as the a posteriori measure on the
environment, after having taken the path $\vv$. 

Define, $\tilde R_N=N^{-1}\sum_{n=1}^N R_n$, with marginals
$\tilde\P_N$ and $\tilde Q_N$. 
Then, since $\VV\times\Omega$ is compact, one can find a subsequence of the 
$\tilde R_N$'s that converges weakly to a probability measure 
$R_\infty$ on $\VV\times\Omega$. In fact, $R_\infty$
will be supported on $\VV_\infty\times\Omega$.

Define now $\P_\infty,Q_\infty$ to be the 
marginals of $R_\infty$ on $\Omega$ and $\VV_\infty$, 
respectively.  Notice that 
\begin{eqnarray*}
\int\!\!P_0^\w(T^{X_1}\w\in A)d\P_n
&=&\int\!\!\sum_{\abs{e}\leq M}\ind(T^e\w\in A)\pi_{0e}(\w)
\sum_{x\in\Z^d}P_x^\w(X_n=0)d\P\cr
&=&\int_A\!\sum_{x\in\Z^d}\sum_{\abs{e}\leq M}\pi_{e0}(\w)P_x^\w(X_n=e)d\P
=\int_A\!\!d\P_{n+1}.
\end{eqnarray*}
This implies that $\P_\infty$ is  an
invariant measure for the process $(T^{X_n}\w)_{n\geq0}$.

Let  $\P_\svv$ be given by the disintegration formula
$$\P_\infty=\int\!\!\P_\svv d Q_\infty(\vv).$$

We would like to show that the conditions of lemma
\ref{invariant_abs2} are in effect. For this, define, for $k\leq0$,
and $\vv\in\cup_{n\geq1}\VV_n$, 
$$A_k(\vv)=\sup_{\w\in\Omega}{\ds{d\P_\svv}_{\big|{\mathfrak{S}_k}}
\over d\P_{\big|{\mathfrak{S}_k}}}(\w).$$
Also, define, for $a>0$, the measure
$$\tilde\theta^{a,k}_N=\int_{A_k\leq a}\P_\svv d\tilde Q_N(\vv).$$
Then, one has that 
${{d\tilde\theta^{a,k}_N}_{|{\mathfrak{S}_k}}\over
d\P_{|{\mathfrak{S}_k}}}\leq a$
and, therefore, one can find a further subsequence of the 
$\tilde\theta^{a,k}_N$'s that converges to a measure 
$\theta^{a,k}_\infty$, with
${{d\theta^{a,k}_\infty}_{|{\mathfrak{S}_k}}\over
d\P_{|{\mathfrak{S}_k}}}\leq a$.
Moreover, one clearly has, for each $N$, 
$\tilde\theta^{a,k}_N\leq\tilde\P_N$. Passing $N$ to
infinity, one has ${\theta^{a,k}_\infty}_{\big|\mathfrak{S}_k}
\leq{\P_\infty}_{\big|\mathfrak{S}_k}$.
Thus, using the same notations as in
lemma \ref{invariant_abs2}, it follows that
\begin{eqnarray}
\label{goal}
\P_\infty^{k,\ll}(\Omega)\geq\theta^{a,k}_\infty(\Omega)
\geq\linf_{N\rightarrow\infty}
\tilde\theta^{a,k}_N(\Omega)=\linf_{N\rightarrow\infty}
\tilde Q_N(A_k\leq a).
\end{eqnarray}
So, according to lemma \ref{invariant_abs2}, to use theorem
\ref{invariant_multi} for the purpose of proving a law of large
numbers, one needs to show that
\begin{eqnarray}
\label{goal2}
\inf_k\sup_a\linf_{N\rightarrow\infty}\tilde Q_N(A_k\leq a)>0.
\end{eqnarray}
Assume now that the mixing condition (\ref{mixing}) holds. 
Then, due to lemma \ref{dobrushin_bound}, one has that,
for $\vv\in\VV_n$,
\begin{eqnarray*}
{\ds{d\P_\svv}_{\big|\mathfrak{S}_k}
\over d\P_{\big|\mathfrak{S}_k}}&=&
\E\left(\left.{\ds\pi_\svv\over\E(\pi_\svv)}\right|\mathfrak{S}_k\right)
\leq\E\left(\left.{\ds\pi_{\svv\cap H_{k-r}^c}\over\E(\pi_\svv)}
\right|\mathfrak{S}_k\right)\cr
&\leq&{\ds\E\left(\pi_{\svv\cap H_{k-r}^c}\right)\over\E(\pi_\svv)}
\exp\left(C\sum_{{}^{x\in\partial_r H_k^c}_{y\in\svv\cap H_{k-r}^c}}
e^{-g\,{\scriptstyle\mathrm{dist}}(x,y)}\right)\cr
&\leq&\kappa^{-{\scriptstyle\mathrm{card}}(\svv\cap H_{k-r})}
\exp\left(\tilde C\sum_{y\in\svv\cap H_{k-r}^c}
e^{-0.5g\,{\scriptstyle\mathrm{dist}}(y,H_k)}\right)\cr
&\leq&\kappa^{-{\scriptstyle\mathrm{card}}(\svv\cap H_{k-r})}
\exp\left(\tilde C\sum_{i\geq r}V_{k-i}(\vv)e^{-0.5gi}\right)=Z_k(\vv),
\end{eqnarray*}
where $\vv\cap H_{k-r}^c=\{x_i\in H_{k-r}^c,0\leq i\leq n\}$, 
$V_j(\vv)=\card{\vv\cap(H_{j-1}\backslash H_j)}$, and
$$\pi_{\svv\cap H_{k-r}^c}=\!\!\!\!\!\!\!\!
\prod_{{\phantom{blabla}}^{i=0}_{x_i\in H_{k-r}^c}}^n
\!\!\!\!\!\!\!\!\!\!\pi_{x_i,x_{i+1}}.$$

Clearly, the left-hand-side in (\ref{goal2}) is bounded from below by
$$\inf_k\sup_a\linf_{N\rightarrow\infty}\tilde Q_N(Z_k\leq a).$$
For a path $(X_n)_{n\geq0}$, define $\tilde Z_{k,n}$ to be
$$\tilde Z_{k,n}=Z_k(X_0-X_n,X_1-X_n,\cdots,X_{n-1}-X_n,0).$$
Also, let $\tau_s=\inf\{n\geq0:X_n.\ell\geq s\}$. Then, for any
$\delta\in(0,1)$, one has
\begin{eqnarray}
\label{lower_estimate}
\tilde Q_N(Z_k\leq a)&=&N^{-1}\sum_{n=1}^N P_0(\tilde Z_{k,n}\leq a)
\geq N^{-1}E_0\left(\sum_{1\leq j\leq\delta N}
\ind_{\tilde Z_{k,\tau_j}\leq a}
\ind_{\tau_{\delta N}\leq N}\right)\cr\cr
&=&N^{-1}\sum_{1\leq j\leq\delta N}P_0(\tilde Z_{k,\tau_j}\leq a)
-N^{-1}\sum_{1\leq j\leq\delta N}
P_0(\tilde Z_{k,\tau_j}\leq a,\tau_{\delta N}>N)\cr\cr
&\geq&N^{-1}\sum_{1\leq j\leq\delta N}P_0(\tilde Z_{k,\tau_j}\leq a)
-\delta P_0(\tau_{\delta N}>N).
\end{eqnarray}

On one hand, we have, 
\begin{eqnarray}
\label{term2}
P_0(\tau_{\delta N}>N)\leq N^{-1}E_0(\tau_{\delta N})\leq
(N\eps)^{-1}E_0(X_{\tau_{\delta N}})\leq{\ds\delta N+M\over N\eps},
\end{eqnarray}
where we have used lemma \ref{kalikow_ballistic}. 
On the other hand, 
$$P_0(\tilde Z_{k,\tau_j}\leq a)
\geq 1-a_1^{-1}E_0(\hat V_{j+k-r,j+M}^j)-
a_2^{-1}\sum_{i\geq r}E_0(\hat V_{j+k-i,j+M}^j)e^{-0.5gi},$$
where $a_1=0.5\log a/\log(\kappa^{-1})$, $a_2=0.5\log a/\tilde C$,
and 
$$\hat V_{i_1i_2}^j=\mathrm{card}\{n:0\leq n\leq\tau_j,i_1\leq X_n.\ell<i_2\}.$$
We had to enlarge the $V_j$'s we had before, to take into account the fact 
that the position of $X_{\tau_j}$ is not known precisely.
To estimate the above expectations, notice that one has, path by path,
$$\sum_{{\phantom{bla}}^{0\leq n\leq\tau_j-1}_{i\leq X_n.\ell<j}}
\!\!\!\!\!\!(X_{n+1}-X_n).\ell\leq(j-i)+M.$$
Using Kalikow's condition (\ref{kalikow}), one has the following
\begin{eqnarray*}
E_0(\hat V_{i,j+M}^j)&=&1+
E_0\left(\sum_{{\phantom{bla}}^{0\leq n\leq\tau_j-1}_{i\leq x.\ell<j}}
\!\!\!\!\!\!\ind_{X_n=x}\right)\leq1+\eps^{-1}
E_0\left(\sum_{{\phantom{bla}}^{0\leq n\leq\tau_j-1}_{i\leq x.\ell<j}}
\!\!\!\!\!\!\ind_{X_n=x}D(T^x\w).\ell\right)\cr
&=&1+\eps^{-1}
E_0\left(\sum_{{\phantom{bla}}^{0\leq n\leq\tau_j-1}_{i\leq X_n.\ell<j}}
\!\!\!\!\!\!(X_{n+1}-X_n).\ell\right)\leq1+\eps^{-1}((j-i)+M).
\end{eqnarray*}
This implies that 
$$P_0(\tilde Z_{k,\tau_j}\leq a)\geq1-a_1^{-1}(1+\eps^{-1}(M+r-k))
-a_2^{-1}\sum_{i\geq r}(1+\eps^{-1}(M+i-k))e^{-0.5gi}.$$
Combining this with (\ref{lower_estimate}) and (\ref{term2}), and taking 
$\delta=0.5\eps$, one has
$$\inf_k\sup_a\linf_{N\rightarrow\infty}
\tilde Q_N(Z_k\leq a)\geq0.25\eps>0.$$
Recalling (\ref{goal}), and using lemma \ref{invariant_abs2}, one has
the existence of the invariant measure $\hat\P_\infty$ that satisfies the
conditions of theorem \ref{invariant_multi}. The transience condition
(\ref{assumption}) is implied by Kalikow's condition (\ref{kalikow}),
due to lemma \ref{infinity}. The law of large numbers, along with the weak convergence of the Cesaro mean of $\P_n$ to $\hat\P_\infty$,
follows from theorem \ref{invariant_multi}.

If the environment is $L$-dependent, instead of mixing, then 
we have 
$${\ds{d\P_\svv}_{\big|\mathfrak{S}_k}
\over d\P_{\big|\mathfrak{S}_k}}\leq
\kappa^{-{\scriptstyle\mathrm{card}}(\svv\cap H_{k-L})},$$
and the rest of the proof is essentially the same as above.

Once one has a law of large numbers, one can use lemma
\ref{kalikow_ballistic}, with $U_L=\{x\in\Z^d:x.\ell\leq L\}$, and Fatou's 
lemma, to show that $T_{U_L}^{-1}X_{T_{U_L}}.\ell\geq LT_{U_L}^{-1}$ 
cannot converge to $0$, proving that the limiting velocity is non-zero.
\qed
\begin{remark}
In the course of preparation of this paper, we learnt of 
\cite{kk}, where the authors prove the law of large numbers for $L$-dependent
non-nestling environments. Their approach is a first step towards the method
we use. They, nevertheless, make use of the regeneration times, introduced in 
\cite{sz}.
Apart from the ellipticity condition, our results include those of \cite{kk}.
We also learnt of \cite{cz}, where the authors use the regeneration
times to prove a result very similar to our theorem \ref{multilln}. However, 
they require moment controls on the regeneration times, which we do not need
in our approach. Working with cones instead of hyperplanes, our method should 
be able to handle mixing on cones, as in \cite{cz}.
\end{remark}

\appendix
\centerline{APPENDIX}\medskip
\refstepcounter{section}
\label{sec_gibbs}
First, we prove a consequence of the Dobrushin-Shlosman mixing
property (\ref{mixing}), in the case of Gibbs fields.
\begin{lemma}
\label{radon_bound_finite}
Let $(\P^\w_V)$ be a Gibbs $r$-specification, corresponding to a
translation-invariant bounded $r$-interaction $U$, and 
satisfying (\ref{mixing}). Then, there exists a constant $C$
such that for all $\Lambda\subset V\subset\Z^d$ 
finite, with $\mbox{dist}(\Lambda,V^c)>r$, 
and for all $x\in V^c$, we have
$$\sup_{\scriptstyle\sigma_{\scriptscriptstyle\Lambda},
\w,\bar\w:\scriptstyle(\w_y)_{y\not=x}=
(\bar\w_y)_{y\not=x}}\abs{{d\P_{V,\Lambda}^\w
\over d\P_{V,\Lambda}^{\bar\w}}(\sigma_\Lambda)-1}\leq 
C\sum_{y\in\partial_r(\Lambda^c)}e^{-g\,{\scriptstyle\mathrm{dist}}(x,y)}.$$
\end{lemma}

\proof Fix $x\in V^c$, and consider $\w,\bar\w\in\Omega$, such that 
$\w_y=\bar\w_y$, for all $y\not=x$. Also, let
$\sigma_\Lambda,\bar\sigma_\Lambda\in\Omega_\Lambda$.
We have, then
$${d\P_{V,\Lambda}^\w \over d\P_{V,\Lambda}^{\bar\w}}
(\sigma_\Lambda)=\E^{\P_V^{\bar\w}}
\left({d\P_V^\w\over d\P_V^{\bar\w}}\bigg|\mathfrak{S}_\Lambda\right)
(\sigma_\Lambda).$$
Notice that, for $\xi_V\in\Omega_V$, we have 
\begin{eqnarray}
\label{radon}
{d\P_V^\w\over d\P_V^{\bar\w}}(\xi_V)=
{\ds\exp\left(-\beta\sum_{A:A\cap V\not=
\emptyset,x\in A}U_A(\w_{V^c},\xi_V)\right)\over\ds\exp\left(-\beta\sum_{A:A\cap V\not=
\emptyset,x\in A}U_A(\bar\w_{V^c},\xi_V)\right)}~.
\end{eqnarray}
So we see that ${d\P_V^\w\over d\P_V^{\bar\w}}$ is 
$\mathfrak{S}_{V_x^r}$-measurable, where 
$V_x^r=\{y\in V:\,\mbox{dist}(x,y)\leq r\}$. Therefore,
$$\E^{\P_V^{\bar\w}}\left({d\P_V^\w\over d\P_V^{\bar\w}}\bigg|
\mathfrak{S}_\Lambda\right)(\sigma_\Lambda)=
\E^{\P_{V-\Lambda,V_x^r}^{\eta}}\left({d\P_V^\w\over 
d\P_V^{\bar\w}}\right),$$
where $\eta=(\bar\w_{\Lambda^c},\sigma_\Lambda)$. Moreover, clearly
\begin{eqnarray}
\label{radon_inf_bound}
\abs{d\P_V^\w\over d\P_V^{\bar\w}}\leq
\exp(2^{{\mbox{\scriptsize{card}}}(\{y\in\Z^d:\norm{y}\leq r\})+1}
\beta\norm{U})=C_1.
\end{eqnarray}
And then, setting 
$\bar\eta=(\bar\w_{\Lambda^c},\bar\sigma_\Lambda)$, we have
\begin{eqnarray}
\label{computation}
\abs{{d\P_{V,\Lambda}^\w\over d\P_{V,\Lambda}^{\bar\w}}
(\sigma_\Lambda)-{d\P_{V,\Lambda}^\w \over d\P_{V,\Lambda}^{\bar\w}}
(\bar\sigma_\Lambda)}&\leq& 
C_1\mbox{Var}(\P_{V-\Lambda,V_x^r}^\eta,
\P_{V-\Lambda,V_x^r}^{\bar\eta})
\leq C_1G\sum_{y\in\partial_r(\Lambda^c)}
e^{-g\,{\scriptstyle\mathrm{dist}}(y,V_x^r)}\cr
&\leq&C_1Ge^{gr}\sum_{y\in\partial_r(\Lambda^c)}
e^{-g\,{\scriptstyle\mathrm{dist}}(y,x)}.
\end{eqnarray}
The conclusion of the lemma follows from integrating out
$\bar\sigma_\Lambda$.\qed\bigskip

Notice now that if (\ref{mixing}) is satisfied, one can define  $\P^\w_V$,
even for an infinite $V$, as the limit of $\P^\w_{V_n}$, for an
increasing sequence of finite volumes $V_n$. It is easy then to see
that (\ref{radon_inf_bound}) still holds, and that, due to the lower
semi-continuity of the variational distance, computation
(\ref{computation}) goes through for all $\Lambda\subset V\subset\Z^d$.
Therefore lemma \ref{radon_bound_finite} still holds for infinite 
$\Lambda\subset V$.\bigskip

{\sc Proof of Lemma} \ref{dobrushin_bound}.
Using $V=H^c$, and applying lemma 
\ref{radon_bound_finite}, we have
$${\E(F|\mathfrak{S}_H)(\w)\over\E(F|\mathfrak{S}_H)(\bar\w)}
={\E^{\P^\w_{V,\Lambda}}(F)\over\E^{\P^{\bar\w}_{V,\Lambda}}(F)}
\leq\!\!\prod_{x\in\partial_r V}
\left(1+C\!\!\!\!\sum_{y\in\partial_r(\Lambda^c)}\!\!
e^{-g\,{\scriptstyle\mathrm{dist}}(x,y)}\right)
\leq\exp\left(C\!\!\!\!\sum_{{}^{x\in\partial_r(H^c)}_{y\in\partial_r(\Lambda^c)}}
e^{-g\,{\scriptstyle\mathrm{dist}}(x,y)}\right).$$
Once again, the conclusion follows by averaging over $\bar\w$.\qed
\bigskip

{\bf Acknowledgments.} 
The author would like to thank his thesis advisor Srinivasa
R. S. Varadhan for suggesting this problem, and for many valuable 
discussions and ideas. The author also thanks O. Zeitouni for 
his comments on earlier versions of this manuscript.
\medskip

\def\refname{\medskip

\hfill
\begin{minipage}{38ex}
{\sc
Courant Institute\hfill\\
New York University\hfill\\
251 Mercer Street\hfill\\
New York, NY, 10003\hfill\\
U.S.A.\hfill\\
E-mail:} {\tt rassoul@cims.nyu.edu}\hfill\\
{\sc URL:} {\tt www.math.nyu.edu/\verb-~-rassoul}\hfill
\end{minipage}
}

\end{document}